\documentclass[12pt]{article}
\usepackage{amsmath}
\usepackage{amsfonts}
\usepackage{amssymb}
\usepackage{theorem}

\title{Ramsey property for spaces with bilinear forms}
\author{Aleksander Ivanov and Fr\'{e}d\'{e}ric Jaffrennou}
\date{ } 

\newtheorem{thm}{Theorem}[section] 
\newtheorem{lem}[thm]{Lemma}

\newtheorem{cor}[thm]{Corollary}
\newtheorem{prop}[thm]{Proposition} 
\newtheorem{remark}[thm]{Remark}

\begin{document}
\maketitle
\topskip 20pt

\begin{quote}
{\bf Abstract.} 
We study the Ramsey property for vector spaces over finite fields with bilinear forms. 
We prove that symplectic spaces over finite fields do not have  the Ramsey property. 
We also describe vector spaces with skew symmetric bilinear forms and radicals of finite codimension, 
where the Ramsey property does not hold. 
Some direct connections with generalized affine spaces are given. 


\end{quote}

\bigskip

%%%%%%%%%%%%%%%%%%%%

\section{Introduction} 
Let $\mathcal{K}$ be a Fra\"{i}ss\'{e} class of finite structures of a finite language and 
$M$ be the Fra\"{i}ss\'{e} limit of  $\mathcal{K}$. 
In particular $\mathcal{K}$ coincides with $Age(M)$, 
the class of isomorphism types of all finite substructures of $M$. 
When $A,B \in \mathcal{K}$ and $A<B$, we denote by ${B\choose A}$ 
the set of all substructures of $B$ which are isomorphic to $A$.  
The class $\mathcal{K}$ is said to have 
the {\bf Ramsey property} if for any $k\in \omega$ 
and a pair $A<B$ from $\mathcal{K}$ 
there exists $C\in \mathcal{K}$ so that 
each $k$-coloring 
$$
\xi :{C\choose A}\rightarrow k
$$ 
is monochromatic on some ${B'\choose A'}$
from $C$ which is a copy of 
${B\choose A}$. 
The latter condition is denoted by  
$$ 
C\rightarrow (B)^{A}_k . 
$$ 
In the situation when $\mathcal{K}$ does not have 
Ramsey property one can consider {\bf Ramsey degrees} of $A$'s  
which are defined as the minimal $k$ such that for every $r\in \omega$ 
and  $B\in \mathcal{K}$ with non-empty ${B\choose A}$
there exists $C\in \mathcal{K}$ so that  
each $r$-coloring 
$$
\xi :{C\choose A}\rightarrow r
$$ 
is $(\le k)$-chromatic on some ${B'\choose A'}$
from $C$ which is a copy of 
${B\choose A}$.
\parskip0pt 

By Theorem 4.8 of the paper of Kechris, Pestov and Todorcevic \cite{KPT},  
if $M$ is a Fra\"{i}ss\'{e} limit of a Fra\"{i}ss\'{e} class $\mathcal{K}$, then the group $Aut(M)$ is 
{\em extremely amenable if and only if the class $\mathcal{K}$ has the Ramsey property and consists of rigid elements.} 
This theorem has started a topic of topological dynamics of automorphism groups of first-order structures, 
which is very active in the last two decades.  
We also mention Theorem 1.2 of A.Zucker from \cite{ZA} which states that  
in the situation above 
{\em conditions 1) and 2) below are equivalent: 

1) $G=Aut(M)$ has mertizable universal minimal flow, 

2) each $A\in \mathcal{K}$ has finite Ramsey degree. } 

\noindent 
It is well known that the age of a countably dimensional vector space over a finite field has the Ramsey property, see \cite{spencer}. 
From this fact one can deduce metrizability of universal minimal flows of 
automorphism groups of $\omega$-categorical $\omega$-stable structures, see \cite{Iv}. 
The argument used in \cite{Iv} suggests that in order to extend this statement to all smoothly approximable structures (\cite{CH})
it is crucial to obtain proerty 2)  of Zucker's theorem above in the case when  $M$ is a vector space with a non-degenerate bilinear form. 
{\em In the present paper we start our investigations of this topic.} 

It is essential that we do not have the restriction of triviality of the radical (as the previous paragraph suggests). 
We have found that spaces with non-trivial radicals (even with radicals of finite codimension) form challenging objects for Ramsey theory. 
This case is studied in Section 2 of this paper. 
It leads to some conjectures directly connected with classical theorems in Ramsey theory of affine/vector spaces.  
{\em We mention Question D in the final part of this section, which looks as a straightforward generalization of the Ramsey property for affine spaces.} 
In Section 1 we prove that symplectic spaces over finite fields do not have the Ramsey property. 
This result answers a question circulated in the community. 
It is obviously connected with the question about Ramsey degrees, formulated in the previous paragraph. 
However in our paper it also motivates the approach presented in Section 2.  

\bigskip 

We now give a short introduction to bilineas spaces. 
From now on we assume that $M$ is a countably infinite dimensional vector space $V$ over a finite field $F$. 
The terminology and notation below are taken from Chapters 19 - 20 of \cite{Asch}. 
Another source is Section 4 of the paper \cite{Barb}. 
The latter paper also discusses connections with the Fra\'{i}ss\'{e} theory. 

Suppose $\sigma \in Aut(F)$. 
A {\bf sesquilinear form} on $V$ is a map 
$\beta: V\times V\rightarrow F$
such that for all $u_i,v_i \in V$, $a,b\in F$
\begin{itemize} 
\item $\beta(u_1 + u_2, v) = \beta(u_1 , v)+ \beta(u_2, v)$, 
\item $\beta(u, v_1 + v_2) = \beta(u, v_1) + \beta(u, v_2)$, 
\item $\beta(au,bv) = a\sigma (b) \beta(u, v)$. 
\end{itemize} 
The form is said to be 
\begin{itemize} 
\item {\bf skew symmetric} if $\sigma = \mathsf{id}$ and 
$\forall u,v\in V  (\beta(u, v)= -\beta(v,u))$, 
\item {\bf symmetric} if $\sigma = \mathsf{id}$ and 
$\forall v\in V (\beta(u, v)= \beta(v,u))$, 
\item {\bf hermitian} if $\sigma \not= \mathsf{id}$, $\sigma^2 = \mathsf{id}$ and $\forall u,v\in V  (\beta(u, v)= \sigma(\beta(v,u))$. 
\end{itemize}  
If $X$ is a subspace of $V$ then define 
\[ 
X^{\perp} = \{ u\in V \, | \, \forall x\in X  (\beta(u, x)=0)\}. 
\]
The {\bf radical} of $V$ (denoted by $\mathsf{Rad}(V)$) is $V^{\perp}$. 
When $\mathsf{Rad}(V) = 0$, the form $\beta$ is called {\bf non-degenerate}. 

When $\beta$ is skew symmetric and non-degenetate, then the pair $(V, \beta)$ is called a {\bf symplectic} space. 
In the case of $char(F) = 2$ we additionally demand that 
$\forall v\in V (\beta(v, v) =0)$. 

A pair of vectors $u,v \in V$ is called a {\bf hyperbolic} pair if $\beta (u,u) = 0=\beta (v,v)$ and $\beta (u,v) = 1$.
It is well-known that any $U<V$ of even dimension and with trivial radical has a basis consisting of pairwise orthogonal hyperbolic pairs. 
\begin{quote} 
{\bf Witt's theorem} (\cite{Asch}, p.81.) 
Let $(V,\beta )$ be a symplectic, orthogonal or unitary space over a finite field $F$ and let $U<V$ be a finite dimensional subspace.  
Then any linear isometry $g:U \rightarrow V$ 
extends to an isometry from $V$ onto $V$.    
\end{quote} 
Since we concentrate on the symplectic case over the two element field, we do not define orthogonal and unitary spaces. 
{\em It is worth noting here that the results of this paper can be adapted to every finite field and to unitary/orthogonal case. 
In order to keep our presentation as easy as possible we usually restrict ourselves by the simplest case.}  
It is worth mentioning that Witt's theorem implies that all these spaces are ultrahomogeneous, i.e.  Fra\'{i}ss\'{e} limits of their ages 
(see \cite{Barb}). 

\begin{remark} \label{ve-sp} 
{\em 
In this remark we mention some consequences of the theorem of Spencer on the Ramsey property for affine spaces. 
Theorem 1 of  \cite{spencer} states that for all $k,r,t \ge 0$ with $t<k$, there exists $n$ with the following property. 
Let $V$ be an $n$-dimensional vector space and $\xi$ be an $r$-coloring of the affine $t$-subspaces of $V$. 
Then there is an affine $k$-space $W\subseteq V$ such that all affine $t$-subspaces of $W$ have the same color.    

Note that there are two types of affine subspaces: vector spaces and spaces of the form $U+v$, 
where $U$ is a vector space and $v\not\in U$. 
We will call the latter form {\em proper affine spaces}. 
The Spencer's statement holds under additional assumptions on the types of the $t$-spaces and of the space $W$.  
} 
The number $n$ above can be chosen as follows. 
Assume that $\xi$ colors proper affine $t$-subspaces of $V$. 
Then $W$ as above can be taken to be a proper affine $k$-space.  
{\em Indeed, applying the original Spencer's theorem to $t$ and $2k$ we find $V$ such that the role of $W$ would play some $W'$ of dimension of $2k$. 
If $W'$ is not proper, then just choose a proper $k$-subspace in it. 

The fact that the theorem holds for vector subspaces is mentioned in \cite{spencer}. 
} 
\end{remark} 

\section{Symplectic spaces over $GF(2)$} 

Let $(V,\beta )$ be symplectic over $GF(2)$ and 
$V= \oplus \{ \langle e_i, e^*_i \rangle \, | \, i\in \omega \}$ be a decomposition of $V$ into hyperbolic subspaces. 
We fix the ordering $e_0 < e^*_0 < \ldots e_i < e^*_i < \ldots$.  
Let $<_{alex}$ be the corresponding anti-lexicographic ordering of $V$ (with respect to the basis $\{ e_i, e^*_i \, | \, i\in \omega\}$): 

$v <_{alex} w$ if comparing elements of the basis which do not appear simultaneously 

both the decomposition of $w$ and of $v$, the former one has a higher element of the 

basis than the elements of the decomposition of $v$ which do not occur in $w$.    \\ 
We consider subspaces of $V$ under the ordering induced by $<_{alex}$. 
The space $V\times V$ (resp. $V\times V \times V$) is considered as a set of pairs under the anti-lexicographic ordering induced by $<_{alex}$ on coordinates.   

Let $U = \langle e_1, e_2 ,e_3 ,e^*_1 + e^*_2 , e^*_3 \rangle$. 
Note that $\mathsf{Rad}(U) = \langle e_1 + e_2 \rangle$. 
Let $W = \langle  e_1, e_2 ,e_3 ,e^*_1 , e^*_2 , e^*_3 \rangle$. 
Note that $\mathsf{Rad}(W) =0$, i.e. $W$ is a symplectic space. 
\begin{remark} 
{\em The map 
\[ 
e_1 \rightarrow e_3\mbox{ , }e_2 \rightarrow e_2 \mbox{ , }  e_3 \rightarrow e_1 \mbox{ , } e^*_1 + e^*_2 \rightarrow e^*_2 + e^*_3 \mbox{ , } e^*_3 \rightarrow e^*_1  
\] 
induces an isometric isomorphism from $U$ onto $\langle e_1 , e_2 , e_3 , e^*_1 , e^*_2 + e^*_3 \rangle$.  
Furthermore, the map 
\[ 
e_1 \rightarrow e_1\mbox{ , }e_2 \rightarrow e_3 \mbox{ , }  e_3 \rightarrow e_2 \mbox{ , } e^*_1 + e^*_2 \rightarrow e^*_1 + e^*_3 \mbox{ , } e^*_3 \rightarrow e^*_2  
\] 
induces an isometric isomorphism from $U$ onto $\langle e_1 , e_2 , e_3 , e^*_2 , e^*_1 + e^*_3 \rangle$. 
By Witt's theorem these isomorphisms are induced by isometric automorphisms of $W$. } 
\end{remark} 

Let us define a coloring 
$c: {(V,\beta ) \choose (U,\beta )} \rightarrow \{ RED , WHITE, BLUE\}$. 
Let $U' <V$ and $(U',\beta ) \cong (U,\beta )$. 
Find in $U'$ the minimal triple $f_1,f_2,f_3$ of pairwise orthogonal elements such that for some $i,j$ the space 
$\langle f_i +f_j \rangle$ coincides with 
$\mathsf{Rad} (U', \beta )$.  
If $\mathsf{Rad} (U', \beta ) = \langle f_1 +f_2 \rangle$
we assign $c(U') = RED$. 
If $\mathsf{Rad} (U', \beta ) = \langle f_1 +f_3 \rangle$
we assign $c(U') = WHITE$. 
Otherwise we color $c(U')= BLUE$. 

\begin{lem}   
If $W'<V$ and $(W',\beta  ) \cong (W,\beta )$, then $W'$ is not monochromatic with respect to $c$. 
\end{lem}

{\em Proof.} 
Find the minimal triple $w_1 ,w_2 , w_3$ of pairwise orthogonal elements of  $W'$. 
Applying Witt's theorem extend the map 
\[ 
e_1 \rightarrow w_1 \mbox{ , } e_2 \rightarrow w_2 \mbox{ , } e_3 \rightarrow w_3 , 
\] 
say $\varphi$, to an isomorphic embedding $\varphi_1 : U \rightarrow W'$ which preserves $\beta$.  
Since $\mathsf{Rad} (\varphi_1 (U))=\langle w_1 + w_2\rangle$, 
then $\varphi_1 (U)$ is $RED$. 
Applying Witt's theorem again, extend $\varphi$ to an isomorphic embedding 
\[ 
\varphi_2 : \langle  e_1, e_2 ,e_3 ,e^*_1 , e^*_2 + e^*_3 \rangle \rightarrow W'. 
\]   
Note that $\langle w_2 + w_3 \rangle$ is the radical of $\varphi_2 (\langle  e_1, e_2 ,e_3 ,e^*_1 , e^*_2 + e^*_3 \rangle)$, i.e. the latter subspace has color $BLUE$.  
$\Box$

\bigskip 

Summarizing the construction and lemma above we obtain the following theorem. 

\begin{thm} 
The age of the symplectic space $(V,\beta )$ does not satisfy the Ramsey property. 
\end{thm}

The following question is open. 

\bigskip 

\noindent 
{\bf Question A.} Do the ages of $\omega$-dimensional symplectic, orthogonal or unitary spaces have finite Ramsey degree (i.e. for every member of the age there is $k$ which bounds the Ramsey degree of it)?

\section{When the number of hyperbolic pairs is bounded}  

Fix a finite number $k\in \omega$. 
Let $(V,\beta )$ be a space over $GF(2)$ with a bilinear skew symmetric form and 
$V= \oplus \{ \langle e_i, e^*_i \rangle \, | \, i\le k \} \oplus \{ \langle e_i  \rangle \, | \, i\in \omega \setminus [ k ]\}$ 
be a decomposition of $V$ into a sum of $k$ hyperbolic subspaces together with an infinite dimensional radical 
$\mathsf{Rad}(V) = \oplus \{ \langle e_i  \rangle \, | \, i\in \omega \setminus [ k ]\}$.  
We denote $V_1= \oplus \{ \langle e_i, e^*_i \rangle \, | \, i\le k \}$. 
The following question is the inspiration of this section. 
\bigskip 

\noindent 
{\bf Question B.} 
Assume that $k\ge 1$. 
Take any $A\in V$ such that the codimension of $A\cap \mathsf{Rad}(V)$  in $A$ is at least 1. 
Does the age of the space $(V,\beta )$ satisfy the Ramsey property with respect to copies of $A$?   

\begin{remark} 
{\em The space $(V,\beta )$ is not an ultrahomogeneous structure: the 1-dimensional space with trivial $\beta$ has copies both in $\mathsf{Rad}(V)$ and in $V_1$. 
On the other hand after naming the elements of $V_1$ as constants, the structure becomes ultrahomogeneous. }  
\end{remark} 
   
\subsection{Possible pairs} 

Let us consider possible pairs $A<B$ of substructures from the age of $V$. 
Let $\pi$ be the projection of $V$ onto  $V_1$. 
It is clear that any finite substructure of $V$ has a 
decomposition of the form $U=U_0 \oplus U_1$ where $U_0 < \mathsf{Rad}(V)$ and  
$U_1$ has the property that $\pi$ is injective on $U_1$. 
Let $A_1 \le  V_1$  be the image of $U_1$ under $\pi$. 
Then every element of $U_1$ is of the form $\mathsf{a} + \mathsf{v}$ with $\mathsf{v}\in \mathsf{Rad}(V)$ and $\mathsf{a}\in A_1$. 
Since the number of possible $A_1$ is finite, every copy of $U$ under an isometry of $V$ belongs to 
the union of finitely many families (not necessarily disjoint) of the following form:  
\[
\mathcal{K}_{C_1} = \{ U_0 \oplus U_1\, | \,  U_0 <\mathsf{Rad}(V)  \mbox{  and }  U_1 \mbox{ is } \pi \mbox{-isometric  with } C{_1} \} ,
\]
\[ 
\hspace{6cm} \mbox{ where } C_1 \le  V_1 \mbox{ is an isometric copy of } A_1 . 
\] 
Note that the set  $\mathcal{K}_{C_1}$ is contained in the the age of $V_{A_1} = \mathsf{Rad}(V) \oplus A_1$. 
It is also worth noting here that when $A_1 =  V_1 \, (= \oplus \{ \langle e_i, e^*_i \rangle \, | \, i\le k \})$, 
every copy of $U$ under an isometry of $V$ 
is a subspace of $V$ which projects onto $A_1$ under $\pi$. 
 
Assuming that $U\in \mathcal{K}_{A_1}$ any coloring 
$\xi :{V\choose U}\rightarrow r$ 
induces some $\xi_{A_1} :{V_{A_1} \choose U}\rightarrow r$. 
This gives plenty pairs $A<B$ demonstrating the failure of the the Ramsey property for the space $(V,\beta )$. 

\begin{lem} \label{l21} 
Let $A< B<V$ be decomposed as $A = A_0 \oplus A_1$ and $B=\mathsf{Rad}(B) \oplus B_1$ with $A_0 < \mathsf{Rad}(V)$ and 
$A_1\le V_1$.    
Assume that $A_1 \le B_1$ and there is an isometric copy of $A_1$ in 
$B_1$ which does not coincide with $A_1$.  
Then there is a coloring $\xi :{V\choose A}\rightarrow r$ which is not monochromatic on any isometric copy of $B$ in $V$.  
\end{lem} 

{\em Proof.} Let $A_1 , A_2 , \ldots A_r$ be the set of  all pairwise distinct copies of $A_1$ in $V_1$. 
For each $i\le r$ we assign color $i$ to all members of $\mathcal{K}_{A_i} \cap {V\choose A}$. 
Then any isometric copy of $B$ has two copies of $A$ of different colors. 
$\Box$ 

\bigskip 

When are the conditions of this lemma satisfied? 
Note that we should avoid the case when $B=\mathsf{Rad}(B)$.  
Indeed, the Ramsey theorem for vector spaces (see \cite{spencer}) works in this situation. 
This means that the only interesting case  is $\mathsf{dim}(B_1 )\ge 2$. 

\begin{thm} \label{boundthm} 
Assume that $k\ge 1$. 
Take any $A, B \in V$ such that $A<B$, $B\not= \mathsf{Rad}(B)$ and $1 \le [A: A\cap \mathsf{Rad}(V) ] < [B: \mathsf{Rad}(B) ]$. 
Then the age of the space $(V,\beta )$ does not satisfy the Ramsey property with respect to copies of $A<B$.  
\end{thm}

\noindent 
{\em Proof.}  
Take $A<B$ as in the formulation and apply Lemma \ref{l21}.  
$\Box$ 

\bigskip 

Consider the case when $B_1$ as above does not have two copies of $A_1$. 
Note that this condition imples that $B_1$ can be embedded into $A_1$. 
Here we have some partial results. 
The first one gives the negative answer under an additional assumption on copies of $B$. 

\begin{prop} \label{l22} 
Assume that $A$ is a finite subspace of $V$ which can be decomposed into $A_0 \oplus A_1$ where $A_1 \le  V_1$  
is non-trivial and $A_0 = A \cap \mathsf{Rad}(V)$. 
Assume that a finite $B\le V_{A_1}$ with $A<B$ has a decomposition $B=\mathsf{Rad}(B) \oplus B_1$ 
with $B_1 \le A_1$ and $\mathsf{dim}(\mathsf{Rad}(B)) \ge 2 \mathsf{dim}(A)$. 
Then there exists a $2$-coloring $\xi :{V\choose A}\rightarrow 2$ which 
is not monochromatic on any ${B'\choose A}$
from $V$ which is a copy of ${B\choose A}$ with $B'_1 \le A_1$.  
\end{prop} 

{\em Proof.} 
Let us fix a base of $A_1$: $a_1, \ldots , a_{\ell}$. 
Take any $U=U_0 \oplus U_1$ where $U_0 < \mathsf{Rad}(V)$ and  
$U_1$ has the property that $\pi$ is injective on $U_1$.  
As we already mentioned above, every element of $U_1$ is of the form $\mathsf{a} + \mathsf{v}$ with $\mathsf{v}\in \mathsf{Rad}(V)$ and $\mathsf{a}\in A_1$. 
Let $u_i = v_i + a_i$, $i\le \ell$, be the base of $U_1$ which corresponds to $a_1 , \ldots , a_{\ell}$ (i.e. $v_i \in V^{\perp}$). 
We color $U$ WHITE if $v_1 , \ldots , v_{\ell}$ are linearly independent over $U_0$. 
Otherwise we color $U$ RED.  
By our assumption on dimensions, it is easy to see that any copy $B'$ of $B$ as in the formulation has both WHITE and RED copies of $A$.  
$\Box$ 

\bigskip

The following observation gives the positive result under additional restrictions on copies of $A$. 

\begin{prop} \label{pRam} 
Assume that $A$ is a finite subspace of $V$ which can be decomposed into $A_0 \oplus A_1$ where 
$A_1 \le  V_1$ and $A_0 = \mathsf{Rad}(A) = A \cap \mathsf{Rad}(V)$. 
Then for any $r$ and any finite $B\le V_{A_1}$ with $A<B$  there exists a finite $C\le V_{A_1}$ so that 
each $r$-coloring $\xi :{C\choose A}\rightarrow r$ 
is monochromatic on the set of all spaces of the form $A'_0 \oplus A_1 \le C$ from some ${B'\choose A}$
which is a copy of ${B\choose A}$ such that $A_1 <B'\le C$.  
\end{prop} 

{\em Proof.} 
Consider $B$ as in the formulation. 
Decompose it as $A_1 \oplus \mathsf{Rad}(B)$. 
Note that $A_0 \le \mathsf{Rad}(B)$.  
Apply the Ramsey theorem for vector spaces from \cite{spencer} to the pair $A_0 \le \mathsf{Rad}(B)$ viewed as a pair from the age of $\mathsf{Rad}(V)$. 
Let $C_0 <\mathsf{Rad}(V)$ be the corresponding subspace given by this theorem. 
Let $C=C_0 \oplus A_1$.  

Take any $r$-coloring of copies of $A$ in $C$ under isometries of $V$. 
This induces an $r$-coloring of  $C_0$-parts of structures decomposed into $A'_0 \oplus A_1$ 
and $A'_0 \subseteq  C_0$.  
By the choice of $C_0$ find $B_2$, a monochromatic copy of $\mathsf{Rad}(B)$ in $C_0$.   
Let $B'= B_2 \oplus A_1$. 
The rest is clear. 
$\Box$

\subsection{Having a distinguished subspace}

The argument of Proposition \ref{pRam} gives the following statement. 

\begin{cor} 
Let $V$ be a vector space over $GF(2)$ of dimension $\omega$ and $A_1$ be a finite subspace of $V$. 
Then the family of all finite subspaces of $V$ which contain $A_1$ has the Ramsey property.  
\end{cor} 

{\em Proof}. Apply the argument of Proposition \ref{pRam} where instead of $\mathsf{Rad}(V)$ (resp. $\mathsf{Rad}(B)$)  we take any complement of $A_1$ in $V$ (resp. the corresponding part in $B$).  $\Box$

\bigskip 

This statement motivates the following dual question. 
\bigskip 

\noindent
{\bf Question C.} 
Let $V$ be a vector space over $GF(2)$ of dimension $\omega$ and $A_1$ be a finite subspace of $V$. 
Let $\pi$ be a projection $V \to A_1$.  
Does the family of all finite subspaces $U<V$ such that $\pi (U)$ contains $A_1$, has the Ramsey property?   
\bigskip 

As above we denote by $\mathcal{K}_{A_1}$ the set of all  finite subspaces $U<V$ such that $\pi (U)$ contains $A_1$. 
Let us consider the case when $\mathsf{dim}(A_1 ) \le 1$. 

\begin{prop} 
Assume that $A_1 <A$ are a finite subspaces of $V$ with $\mathsf{dim}(A_1 )\le 1$. 
Then for any $r$ and any finite $B\in \mathcal{K}_{A_1}$ with $A <B$  there exists a finite $C\in \mathcal{K}_{A_1}$ so that 
each $r$-coloring $\xi :{C\choose A}\rightarrow r$ 
is monochromatic on the subset of  $\mathcal{K}_{A_1}$ of the form ${B'\choose A}$
from $C$ which is a copy of ${B\choose A}$.  
\end{prop} 

{\em Proof.} 
We start with the case $\mathsf{dim}(A_1 ) =1$.  
Let $V= V_0 \oplus A_1$ where $V_0 = \mathsf{ker} (\pi )$. 
Any finite subspace from $\mathcal{K}_{A_1}$ has a decomposition of the form $U=U_0 \oplus U_1$ where 
$U_0 < \mathsf{ker}(\pi )$ and  
$U_1$ has the property that $\pi$ is injective on $U_1$ and $\pi (U_1 ) = A_1$. 
Let us fix a base of $A_1$, say $a_1$.   
Then a generator of $U_1$ is of the form $a_1 + v_1$ with some $v_1 \in V_0$. 
This induces an assignment 
\[
U \rightarrow (U_0 , v_1  )  \to U_0 + v_1 . 
\] 
Note that  $U_0 + v_1$ is uniquely determined by $U$ and $a_1$: 
changing $U_1$ in the decomposition, the assignment above does not change the space $U_0 + v_1$. 
Choose a copy of $A$ of the form $(V_0 \cap A) \oplus \langle v_1 +a_1 \rangle$ where $v_1 \not\in V_0 \cap B$. 
Now apply the Ramsey theorem for affine spaces from \cite{spencer} to copies of $(V_0 \cap A) +v_1$ 
(i.e. proper affine spaces of dimension $\mathsf{dim} (V_0 \cap A)$, see Remark \ref{ve-sp}) 
and their embeddings into copies of the corresponding $(B\cap V_0 )+v_1$. 
We obtain an affine space $C_0 + v'$ such that for any coloring of copies of $(V_0 \cap A)+v_1$, 
there are monochromatic copies of $(B \cap V_0 )+ v_1$ in it.  
Let $C = C_0 \oplus \langle a_1 +v'\rangle$. 
Then $C\in \mathcal{K}_{A_1}$. 

Note that for any $\xi$ as in the formulation of the proposition, and $U<C$, a copy of $A$, 
we have a well defined color of the affine space $U_0 + v_1$ 
computed as above: just take the color $\xi (U )$. 
As a result we have a coloring of affine subspaces of $C_{0}+v'$ of dimension $\mathsf{dim}(V_0 \cap A)$ which are not vector spaces.  
Choosing a monochromatic copy of $(B\cap V_0 )+v_1$, say $W + v'_1$, we obtain a monochromatic $W \oplus\langle a_1 + v'_1 \rangle$ in $C$.  

The case when $\mathsf{dim}(A_1 ) = 0$ is similar. 
In this case we just consider subspaces $U< V$. 
The Ramsey theorem for vector spaces from \cite{spencer} gives the result. 
$\Box$

\bigskip 

The argument which is applied in the proof of this proposition leads to the following problem. 

\bigskip 

\noindent 
{\bf Question D.}  
Let $V$ be a vector space of dimension $\omega$ over a finite field. 
Let $n\in \omega$. 
When $U$ is a finite subspace of $V$ and $v_1 , \ldots , v_n$ is a sequence linearly independent over $U$ we call the tuple 
$(U+v_1 , \ldots , U+v_n )$ an $n$-{\em space-tuple based on} $U$.  
When $(W+v'_1 , \ldots , W+ v'_n )$ is another $n$-space-tuple, then by  $(U+v_1 , \ldots , U+v_n ) < (W+v'_1 , \ldots , W+ v'_n )$ we denote the situation when each $U+v_i$ is a subspace of $W+v'_i$. 
{\em Does the family of $n$-space-tuples of the space $V$ have the Ramsey property with respect to substructure relation $\le$?} 
Here we exactly mean the following statement. 
For all $k,r,t \ge 0$ with $t<k$, there exists $m$ with the following property. 
Let $V$ be an $m$-dimensional vector space and $\xi$ be an $r$-coloring of the $n$-space-tuples based on $t$-subspaces of $V$. 
Then there is an $n$-space-tuple in $V$ of dimension $k$, say $(W +v_1 , \ldots, W+v_n)$, such that all its $n$-space subtuples of dimension $t$ have the same color.    
\bigskip 

Note that in the case $n=1$ this is just the Ramsey theorem for affine spaces. 
We now show that the statement of this problem implies the positive solution to the question formulated in the beginning of this subsection, i.e. Question C. 

\begin{prop} 
Assume that $A_1 <A$ are a finite subspaces of $V$ with $\mathsf{dim}(A_1 )=n$. 
If the family of $n$-space-tuples of an infinitely dimensional space has the Ramsey property, 
then for any $r$ and any finite $B\le \mathcal{K}_{A_1}$ with $A <B$  there exists a finite $C\in \mathcal{K}_{A_1}$ so that 
each $r$-coloring $\xi :{C\choose A}\rightarrow r$ 
is monochromatic on the subset of  $\mathcal{K}_{A_1}$ of the form ${B'\choose A}$
from $C$ which is a copy of ${B\choose A}$. 
\end{prop} 

{\em Proof.} 
We apply the argument already used in the case $\mathsf{dim}(A_1 ) =1$.  
Let $V= V_0 \oplus A_1$, where $V_0 = \mathsf{ker} (\pi )$. 
Any finite subspace of $\mathcal{K}_{A_1}$ has a decomposition of the form $U=U_0 \oplus U_1$ where 
$U_0 < V_0$ and  
$U_1$ has the property that $\pi$ is injective on $U_1$ and $\pi (U_1 ) = A_1$. 
Let fix a base of $A_1$, say $a_1 , \ldots , a_n$.   
Then there is a base of $U_1$ of the form $a_1 + v_1 , \ldots , a_n +v_n$ with some 
$v_1, \ldots, v_n  \in V_0$. 
This induces an assignment 
$$
U \rightarrow (U_0 +v_1 , \ldots U_0 + v_n ). 
$$  
(The resulting tuple is not necessary an $n$-space-tuple.)  
This map is determined by $U_0$ and $a_1 ,\ldots , a_n$, and is injective. 
Indeed, any decomposition of $U$ as above is of the form 
$$
U_0 \oplus  \langle a_1 + v_1 +u_1 , \ldots , a_n +v_n +u_n \rangle , 
$$ 
where $u_1 , \ldots , u_n\in U_0$.   
Thus  $(U_0 +v_1 , \ldots U_0 + v_n )$ is defined uniquely. 

Choose a copy of $A$ of the form $(V_0 \cap A) \oplus \langle v_1 +a_1 , \ldots , v_n + a_n \rangle$, 
where $v_1 , \ldots v_n $ are linearly independent over $V_0 \cap B$. 
Now apply our assumption that the Ramsey theorem for $n$-space-tuples holds. 
We consider colorings of  copies of the $n$-space-tuple $((V_0 \cap A) +v_1 , \ldots ,(V_0 \cap A)+v_n )$ 
together with their embeddings into copies of the corresponding $((B\cap V_0 )+v_1 , \ldots ,(B\cap V_0 )+v_n )$. 
We obtain a space $C_0$ such that for any coloring of $n$-space-tuples from it of dimension $\mathsf{dim}(V_0 \cap A)$ 
there are monochromatic copies of $(B_0 + v_1 , \ldots , B_0 +v_n )$ in it. 
We may assume that $\bigcup \{ (B\cap V_0 ) + v_i | i\le n \}  \subseteq  C_0$.  
Let $C = C_0 + \langle a_1 , \ldots , a_n  \rangle$. 

Note that for any $\xi$ as in the formulation of the proposition, and $U<C$, a copy of $A$, 
we have a well defined color of the $n$-space-tuple $(U_0 + v_1 , \ldots , U_0 +v_n )$ 
computed as above:  
$$ 
(U_0 + v_1 , \ldots , U_0 +v_n ) \rightarrow \xi (U ) .  
$$
As a result we have a coloring of $n$-space-tuples from $C_{0}$ of dimension $\mathsf{dim}(V_0 \cap A)$.  
Choosing a monochromatic $n$-space-tuple of dimension $\mathsf{dim}(B\cap V_0 )$, say $(W + v'_1 .\ldots , W+v'_n )$, 
we obtain a monochromatic $W \oplus\langle a_1 + v'_1 ,\ldots , a_n + v'_n \rangle$ in $C$.  
$\Box$

\vspace*{10mm}

\begin{flushleft}
\begin{footnotesize}
Aleksander Iwanow, \,  
University of Opole, \, 
Institute of Computer Science, \\
ul. Oleska 48, \,  45 - 052,  Opole, Poland \\ 
email: aleksander.iwanow@uni.opole.pl
\end{footnotesize} 

\bigskip 

\begin{footnotesize}
Fr\'{e}d\'{e}ric Jaffrennou, \,
Pozna\'{n}  University of Technology, \,
Institute of Mathematics, \\ 
 ul. Piotrowo 3A, \,  61--138 , Pozna\'{n}, \, 
Poland\\
email: frederic.jaffrennou@put.poznan.pl
\end{footnotesize}

\end{flushleft}

\end{document}